\documentclass[12pt,reqno]{amsart}

\usepackage{tikz}
\usepackage{times}
\usepackage[T1]{fontenc}
\usepackage{mathrsfs}
\usepackage{latexsym}
\usepackage{hyperref}

\newtheorem{thm}{Theorem}[section]

\newtheorem{rek}[thm]{Remark}

\numberwithin{equation}{section}

\colorlet{shadecolor}{gray!30}
\usepackage{breakurl}

\begin{document}

\title{Squaring the Circle Revisited}

\author{H\`ung Vi\d{\^e}t Chu}

\email{\textcolor{blue}{\href{mailto:hungchu2@illinois.edu}{hungchu2@illinois.edu}}}
\address{Department of Mathematics\\ University of Illinois Urbana-Champaign, Urbana, IL 61820, USA}


\begin{abstract}
    Squaring the circle is impossible, but it can be squared approximately. Ramanujan gave a construction correct to eight decimal places. In his book \textit{Mathographics}, Dixon gave constructions correct to three decimal places. In this article, we provide a new construction correct to three decimal places and another correct to nine decimal places. 
\end{abstract}

\maketitle

\section{History}
Ancient geometers asked whether it was possible to construct a square that has the same area as a given circle using only straightedges and compasses. Because it is possible to approximate the area of a circle with inscribed regular polygons to any degree of precision, and polygons can be squared, it is reasonable to expect that a circle can be squared as well. While many Greek mathematicians and philosophers including Anaxagoras, Hippocrates of Chios, and Oenopides attempted a solution, mathematicians after the sixteenth century tried to prove its impossibility. However, it was not until the nineteenth century that considerable progress was made. In 1837, Pierre Wantzel \cite{Wa} proved that lengths that can be constructed must be zeros of polynomials with integer coefficients. This was a breakthrough towards proving the impossibility of a solution; that is, it suffices to prove that $\pi$ is not a solution of a polynomial with integer coefficients (or transcendental over $\mathbb{Z}$, in other words). Almost 50 years later, Lindemann \cite{Lin} proved the transcendence of $\pi$ and thus, closed the problem.

However, mathematicians continued to square the circle approximately. As we will discuss soon, the problem of constructing a square of area $\pi$ is the same as constructing either $\pi$ or $\sqrt{\pi}$. Ramanujan \cite{Ra} constructed $\frac{355}{113}$, an approximation of $\pi$ correct to six decimal places. Later, he \cite{Ra2} gave another construction correct to eight decimal places, where he used $\sqrt[4]{9^2+19^2/22}$ to approximate $\pi$. An image of his construction is available at: \url{https://commons.wikimedia.org/wiki/File:01-Squaring_the_circle-Ramanujan-1914.gif}. Because there are rational numbers arbitrarily close to $\pi$, we can construct lengths approximating $\pi$ up to any degree of precision. Unfortunately, there is a tradeoff between the level of precision and the complexity of the construction. Higher precision often requires more steps involved; the best construction should balance these two factors. In 1991, Dixon gave constructions for $6/5\cdot (1+\phi)\approx 3.1416$, where $\phi = \frac{1+\sqrt{5}}{2}$, the golden ratio \cite{D}. In this article, we show a new method of constructing $6/5\cdot (1+\phi)$ and then construct $\sqrt{\frac{63}{25}(1+\frac{5}{2}\frac{15\sqrt{5}-7}{269})}$ (a number due to Ramanujan \cite{Ra2}) whose square approximates $\pi$ correct to nine decimal places. 

Before moving on, we review some basic straightedge-and-compass constructions. The following constructions are numbered for later reference.
\begin{enumerate}
    \item \label{perline}Given a line $\ell$ and a point, it is possible to draw a line that goes through the point and is perpendicular to $\ell$.
    \item \label{div} Let $n\in \mathbb{N}$. Given a line segment of length $a$, it is possible to divide it into $n$ equal segments. Therefore, it is possible to construct a new segment with length equal $ra$ for all $r\in \mathbb{Q}$. 
    \item Given a line segment, it is possible to draw a circle with it as a diameter. 
    \item \label{prod} Given two line segments with lengths $a$ and $b$, we are able to construct segments of lengths $a\cdot b$ and $a/b$.
    \item \label{sqroot} Given a line segment of length $a$, we are able to construct a segment of length $\sqrt{a}$.
\end{enumerate}
All of these basic constructions can be found in \textit{Mathographics} by Dixon \cite{D}. It is worth noting that due to \eqref{perline}, \eqref{prod}, and \eqref{sqroot}, the problem of squaring the circle is the same as constructing either $\pi$ or $\sqrt{\pi}$.

\section{Squaring the Circle Using $6/5\cdot (1+\phi)$}
A naive way to construct $\sqrt{6/5(1+\phi)}$ would be to construct $\phi$, then add $1$ to it, then multiply the whole length by $6/5$, and finally take the square root of the resulting length. However, such a construction is not elegant. Figure \ref{ref:4dp} shows our construction. Here are the steps:
\begin{enumerate}
    \item[(i)] Let $AB = 1$. Construct $BC \perp AB$ and $BC = 2$.
    \item[(ii)] Draw the circle centered at $A$, radius $AC$, which cuts $AB$ extended at $D$. 
    \item[(iii)] Let $M$ be on $AD$ such that $DM = 2AD/5$. 
    \item[(iv)] Let $N$ be such that $ND = DM/2$.
    \item[(v)] Draw the circle taking $NB$ to be one of its diameters. 
    \item[(vi)] Let $H$ be on the circle such that $MH \perp NB$. Then $MH = \sqrt{6/5\cdot (1+\phi)}$, our desired quantity. 
\end{enumerate}
\begin{figure}[htb]
\centering
\begin{tikzpicture} 
[acteur/.style={circle, fill=black,thick, inner sep=1pt, minimum size=0.1cm}] 
\draw (1.6,0)[black, fill = shadecolor] circle (1.6 cm);
\draw (-2.1462,0)[black, fill = shadecolor] rectangle (-4.9817,-2.8355);
\node (a1) at (0,0) [acteur][label = below: A]{};
\node (a2) at (1.6,0) [acteur][label = right: B]{};
\node (a3) at (1.6,3.2) [acteur][label = right: C]{};
\node (a4) at (-3.577,0) [acteur][label = below: D]{};
\node (a5) at (-2.1462,0) [acteur][label = below: M]{};
\node (a6) at (-4.2924,0) [acteur][label = below: N]{};
\node (a7) at (-2.1462, -2.8355) [acteur][label = below: H]{};
\draw (a1)--(a3);
\draw (a6)--(a4);
\draw (a2)--(a3);
\draw (a4)--(a5)--(a1)--(a2);
\draw (a1)[black]  circle (3.5777 cm);
\draw (-1.3462,0)[black] circle (2.9462 cm);
\draw (a5)--(a7);
\end{tikzpicture}
\caption{The shaded square and the shaded circle almost have the same area.}
\label{ref:4dp}
\end{figure}

We now prove that $MH = \sqrt{6/5\cdot (1+\phi)}$.

\begin{proof}
Because $D$ and $C$ are on the same circle centered at $A$, $AD = AC = \sqrt{1^2+2^2} = \sqrt{5}$. Because $DM = 2AD/5, DM = 2/\sqrt{5}$ and $MA = 3/\sqrt{5}$. So, $MN = 3/\sqrt{5}$, and $MB = 3/\sqrt{5}+1$. Since $NB$ is the diameter of a circle and $H$ is a point on the circle, $\triangle NHB$ has a right angle at $H$. Therefore, 
$$MH \ =\ \sqrt{NM\cdot MB} \ =\ \sqrt{3/\sqrt{5}\cdot (3/\sqrt{5}+1)} \ =\ \sqrt{6/5\cdot (1+\phi)},$$
as desired. 
\end{proof}
\begin{rek}\normalfont
We compute the error
$$\frac{6/5\cdot (1+\phi)}{\pi} \ \approx\ 1.000015.$$
So, for every one million parts, we are off by about fifteen parts. This approximation is very good!
\end{rek}

\section{Squaring the Circle Using $\frac{63}{25}\big(1+\frac{5}{2}\cdot \frac{15\sqrt{5}-7}{269}\big)$}
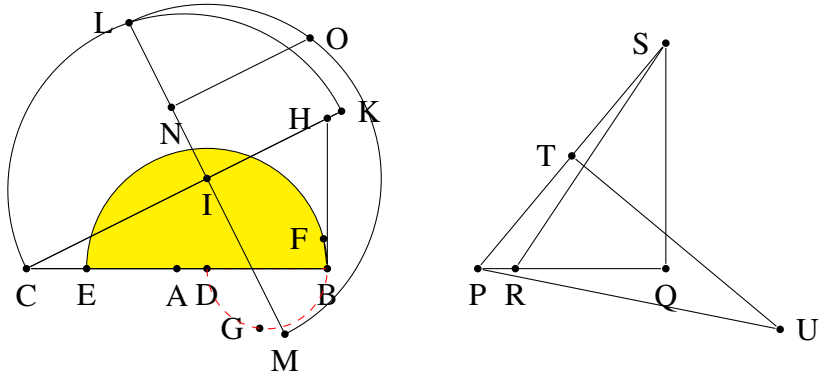
\begin{figure}[htb]
\centering
\begin{tikzpicture}
[acteur/.style={circle, fill=black,thick, inner sep=1pt, minimum size=0.1cm}] 
\draw [rotate = 0,  fill = shadecolor] (2,0) arc(0:180:1.6 ) -- cycle;
\node (a1) at (0,0) [acteur][label = below: A]{};
\node (a2) at (2,0) [acteur][label = below: B]{};
\node (a3) at (-2,0) [acteur][label = below: C]{};
\node (a4) at (0.4, 0) [acteur][label = below: D]{};
\node (a5) at (-1.2, 0) [acteur][label = below: E]{};
\node (a6) at (1.95,0.3968) [acteur][label = left: F]{};
\node (a7) at (1.1,-0.7937) [acteur][label = left: G]{};
\node (a8) at (2,2) [acteur][label = left: H]{};
\node (a9) at (0.4,1.2)  [acteur][label = below: I]{};
\node (a10) at (2.1888,2.0944) [acteur][label = right: K]{};
\node (a11) at (-0.635992, 3.27198) [acteur][label = left: L]{};
\node (a12) at (1.435992,-0.87198) [acteur][label = below: M]{};
\node (a13) at (-0.07328,2.1466)  [acteur][label = below: N]{};
\node (a14) at (1.76977,3.068125) [acteur][label = right: O]{};
\node (a15) at (4,0) [acteur][label = below: P]{};
\node (a16) at (6.5,0) [acteur][label = below: Q]{};
\node (a17) at (4.5,0) [acteur][label = below: R]{};
\node (a18) at (6.5,3)  [acteur][label = left: S]{};
\node (a19) at (5.25,1.5) [acteur][label = left: T]{};
\node (a20) at (8.02,-0.80822) [acteur][label = right: U]{};

\draw (a19)--(a20);
\draw (a15)--(a20);
\draw (a14) -- (a13);
\draw (a1)--(a2);
\draw (a1)--(a3);
\draw [rotate = 180, red, dashed] (a4) arc(0:180:0.8) -- cycle;
\draw (a6)--(a2);
\draw (a8)--(a2);
\draw (a8)--(a3);
\draw (a10)--(a8);
\draw [rotate = 26.565] (a10) arc(0:180:2.34161 ) -- cycle;
\draw[rotate = -63.435] (a12) arc(0:180:2.3165) -- cycle;
\draw (a15)--(a16);
\draw (a15)--(a18);
\draw (a17)--(a18);
\draw (a16)--(a18);
\draw (a2)--(a3);
\end{tikzpicture}
\caption{Construction of $\frac{\sqrt{63}}{5}$, $\frac{\sqrt{15\sqrt{5}-7}}{5}$, and $\frac{\sqrt{269}}{8}$}
\label{ref:9dp}
\end{figure}

Note that 
\begin{align*}
&\frac{63}{25}\bigg(1+\frac{5}{2}\cdot \frac{15\sqrt{5}-7}{269}\bigg) \ \approx\ 3.1415926538,\\
\mbox{while } &\pi \ \approx\ 3.1415926535.
\end{align*}
This approximation is correct to nine decimal places.  Of course, a naive way would be to construct a length of $\sqrt{5}$ and then lengthen it by $15$ times to then divide it by $269$. We show a more elegant construction. Figure \ref{ref:9dp} shows how we construct $\frac{\sqrt{63}}{5}$, $\frac{\sqrt{15\sqrt{5}-7}}{5}$, and $\frac{\sqrt{269}}{8}$. Once we have these lengths, we can construct $\big(\frac{63}{25}\big(1+\frac{5}{2}\cdot \frac{15\sqrt{5}-7}{269}\big)\big)^{1/2}$ with basic constructions. 

Here are the steps:

Left figure: 
\begin{enumerate}
    \item[(i)] Let $AB = AC = 1$. Let $D$ be on $AB$ so that $AD = AB/5$. 
    \item[(ii)] Draw a semicircle with radius $DB$ that intersects $BC$ at $E$. Pick $F$ on the semicircle such that $BF = AD$.
    \item[(iii)] Draw a semicircle with diameter $DB$. Pick $G$ on the semicircle so that $BG = EA$.
    \item[(iv)] Let $H$ be such that $HB \perp  BC$ and $HB = 1$.
    \item[(v)] Let $I$ be on $CH$ such that $CI = 3CH/5$. On $CH$ extended, pick $K$ so that $IK = 1$.
    \item[(vi)] Draw a semicircle with diameter $CK$. Let $L$ be on the semicircle so that $IL \perp CK$. Pick $M$ so that $I$ is the midpoint of $LM$.
    \item[(vii)] Pick $N$ on $IL$ so that $NI = DG$.
    \item[(viii)] Draw a semicircle with diameter $LM$. Pick $O$ on the semicircle so that $NO \perp LM$. 
\end{enumerate}

Right figure: 
\begin{enumerate}
    \item[(i)] Let $RQ = 1$ and $PR = 1/4$.
    \item[(ii)] Let $SQ = 3/2$ and $SQ \perp RQ$. 
    \item[(iii)] Let $T$ be the midpoint of $SP$. Draw $UT = SR$ and $UT \perp PS$. Connect $PU$. 
\end{enumerate}

We will prove that $EF = \frac{\sqrt{63}}{5}$, $ON = \frac{\sqrt{15\sqrt{5}-7}}{5}$, and $PU = \frac{\sqrt{269}}{8}$. But for now, assume that they are true. By constructions \eqref{div} and \eqref{prod}, we can construct a length of $\frac{5}{8}\cdot \frac{ON}{PU} = \sqrt{\frac{15\sqrt{5}-7}{269}}$. By these same constructions, we can construct $\frac{5}{2}\cdot \frac{15\sqrt{5}-7}{269}+1$. Then, by construction \eqref{sqroot}, we construct $\sqrt{\frac{5}{2}\cdot \frac{15\sqrt{5}-7}{269}+1}$. Finally, by construction \eqref{prod}, we construct $\frac{\sqrt{63}}{5}\sqrt{\frac{5}{2}\cdot \frac{15\sqrt{5}-7}{269}+1}$.

\begin{proof}
We have
$$EF^2 \ =\ EB^2 - BF^2 \ =\ 4 DB^2-AD^2 \ =\ 4\cdot (4/5)^2 - (1/5)^2 \ =\ 63/25.$$
So, $EF = \frac{\sqrt{63}}{5}$. Since $GB = EA = ED - AD = 4/5-1/5 = 3/5$, we know that $$DG \ =\ \sqrt{DB^2-GB^2} \ =\ \sqrt{(4/5)^2-(3/5)^2} \ =\ \sqrt{7}/5.$$ Also, 
\begin{align*}
    CH \ =\ \sqrt{BC^2+BH^2} \ =\ \sqrt{2^2+1^2} \ =\ \sqrt{5}
\end{align*}
implies $CI = 3CH/5 = 3/\sqrt{5}$. Because $\triangle LCK$ has a right angle at $L$, $LI = \sqrt{CI\cdot IK} =\sqrt{3}/\sqrt[4]{5}$. So, 
\begin{align*}
    LN &\ =\ LI - NI \ =\ LI - DG \ =\ \sqrt{3}/\sqrt[4]{5}-\sqrt{7}/5\\
    NM &\ =\ NI + IM \ =\ LI+DG \ =\  \sqrt{3}/\sqrt[4]{5} + \sqrt{7}/5.
\end{align*}
Because $\triangle LOM$ has a right angle at $O$, $$NO = \sqrt{LN\cdot NM} \ =\ \sqrt{\big(\sqrt{3}/\sqrt[4]{5}-\sqrt{7}/5\big)(\sqrt{3}/\sqrt[4]{5} + \sqrt{7}/5)}\ =\ \frac{\sqrt{15\sqrt{5}-7}}{5}.$$
Lastly, we prove that $PU = \frac{\sqrt{269}}{8}$. We have 
\begin{align*}
    RS^2 &\ =\ RQ^2 + QS^2 \ =\ 1^2 + (3/2)^2 \  =\ 13/4,\\
    PS^2 &\ =\ PQ^2+SQ^2 \ =\ 61/16.
\end{align*}
Hence, $TP^2 = PS^2/4 = 61/64$. Therefore, $$PU^2 \ =\ TP^2+TU^2 \ =\ TP^2+RS^2 \ =\ 13/4+61/64 \ =\ 269/64$$ and so $PU = \frac{\sqrt{269}}{8}$. 
\end{proof}

\begin{rek}
\normalfont We compute the error
$$\frac{\frac{63}{25}\big(1+\frac{5}{2}\cdot \frac{15\sqrt{5}-7}{269}\big)}{\pi} \ \approx\ 1.000000000068.$$
So, for every one trillion $(10^{12})$ parts, we are off by about 68 parts. This approximation is extremely good!
\end{rek}

\end{document}